\documentclass[12pt]{amsart}
\usepackage{amsmath,amssymb,amsthm}
\usepackage[]{latexsym}
\newtheorem{defn0}{Definition}[section]
\newtheorem{prop0}[defn0]{Proposition}
\newtheorem{thm0}[defn0]{Theorem}
\newtheorem{lemma0}[defn0]{Lemma}
\newtheorem{corollary0}[defn0]{Corollary}
\newtheorem{example0}[defn0]{Example}
\newtheorem{remark0}[defn0]{Remark}
\newtheorem{conjecture0}[defn0]{Conjecture}
\newtheorem{notation0}[defn0]{Notation}

\newenvironment{definition}{\begin{defn0}\rm}{\end{defn0}}
\newenvironment{proposition}{\begin{prop0}}{\end{prop0}}
\newenvironment{theorem}{\begin{thm0}}{\end{thm0}}
\newenvironment{lemma}{\begin{lemma0}}{\end{lemma0}}
\newenvironment{corollary}{\begin{corollary0}}{\end{corollary0}}

\newenvironment{remark}{\begin{remark0}\rm}{\end{remark0}}
\newenvironment{conjecture}{\begin{conjecture0}}{\end{conjecture0}}

\newcommand{\cN}{{\mathcal N}}
\newcommand{\cS}{{\mathcal S}}
\newcommand{\cD}{{\mathcal D}}

\newcommand{\Lie}{{\mathrm {Lie}}}

\newcommand{\cI}{{\mathcal I}}

\newcommand{\disc}{{\mathrm {disc }}}
\newcommand{\sk}{\vspace{0.1in}}

\newcommand{\Pic}{\mathrm{Pic}}

\newcommand{\Norm}{\mathrm{N}}
\newcommand{\M}{\mathrm{M}}
\newcommand{\Aut}{\mathrm{Aut}}

\newcommand{\GL}{{\mathrm{GL}}}
\newcommand{\End}{{\mathrm{End}}}

\newcommand{\Z}{{\mathbb Z}}
\newcommand{\Q}{{\mathbb Q}}
\newcommand{\C}{{\mathbb C}}
\newcommand{\R}{{\mathbb R}}
\newcommand{\F}{{\mathbb F}}

\newcommand{\cL}{{\mathcal L}}

\newcommand{\cP}{{\mathcal P}}
\newcommand{\cO}{{\mathcal O}}
\newcommand{\NS}{{\mathrm{NS}}}
\newcommand{\n}{{\mathrm{n}}}
\newcommand{\tr}{{\mathrm{tr}}}

\newcommand{\ra}{{\rightarrow}}

\include{thebibliography}

\begin{document}

\footnote{Partially supported by a grant FPI from Ministerio de
Ciencia y Tecnolog\'{\i}a, by MCYT BFM2000-0627 and by DGCYT
PB97-0893. }

\address{Dpt.\,\`{A}lgebra i Geometria, Universitat de Barcelona,
Gran Via, 585, E-08007, Barcelona.}

\email{vrotger@mat.ub.es}

\subjclass{11G18, 14G35}

\keywords{Abelian variety, polarization, quaternion algebra}

$$
\mbox{\bfseries {{\Large Quaternions, polarizations and class
numbers}}}
$$

$$
\mbox{By {\em V. Rotger} at Barcelona}
$$
$\\ $ {\bfseries Abstract.}\quad We study abelian varieties $A$
with multiplication by a totally indefinite quaternion algebra
over a totally real number field and give a criterion for the
existence of principal polarizations on them in pure arithmetic
terms. Moreover, we give an expression for the number $\pi _0(A)$
of isomorphism classes of principal polarizations on $A$ in terms
of relative class numbers of CM fields by means of Eichler's
theory of optimal embeddings. As a consequence, we exhibit simple
abelian varieties of any even dimension admitting arbitrarily many
non-isomorphic principal polarizations. On the other hand, we
prove that $\pi _0(A)$ is uniformly bounded for simple abelian
varieties of odd square-free dimension.

\section{{\bfseries Introduction}}
\sk\noindent

It is well-known that every elliptic curve $E$ over an arbitrary
algebraically closed field admits a unique principal polarization
up to translations. This is in general no longer shared by higher
dimensional abelian varieties and it is a delicate question to
decide whether a given abelian variety $A$ is principally
polarizable. Even, if this is the case, it is an interesting
problem to investigate the set $\Pi _0(A)$ of isomorphism classes
of principal polarizations on $A$. By a theorem of Narasimhan and
Nori (cf.\,\cite{NaNo}), $\Pi _0(A)$ is a finite set. We shall
denote its cardinality by $\pi _0(A)$.

The aim of this paper is to study these questions on abelian
varieties with quaternionic multiplication. It will be made
apparent that the geometric properties of these abelian varieties
are encoded in the arithmetic of their ring of endomorphisms. The
results of this paper shed some light on the geometry and
arithmetic of the Shimura varieties that occur as moduli spaces of
abelian varieties with quaternionic multiplication and their
groups of automorphisms. In this regard, we refer the reader to
\cite{Ro1} and \cite{Ro2}. Our results are also the basis of a
study of the diophantine properties of abelian surfaces with
quaternionic multiplication over number fields carried by
Dieulefait and the author in \cite{DiRo}.

Let us remark that a generic principally polarizable abelian
variety admits a single class of principal polarizations. In
(\cite{Hu}), Humbert was the first to exhibit simple complex
abelian surfaces with two non-isomorphic principal polarizations
on them. Later, Hayashida and Nishi (\cite{HaNi} and \cite{Ha})
computed $\pi _0(E_1\times E_2)$ for isogenous elliptic curves
$E_1/\C $ and $E_2/\C $ with complex multiplication. In positive
characteristic, Ibukiyama, Katsura and Oort (\cite{IbKaOo})
related the number of principal polarizations on the power $E^n$
of a supersingular elliptic curve to the class number of certain
hermitian forms. With similar methods, Lange (\cite{La}) produced
examples of simple abelian varieties of high dimension with
several principal polarizations. However, he showed that for an
abelian variety with endomorphism algebra $\End (A)\otimes \Q =F$,
a totally real number field, the number $\pi _0(A)$ is uniformly
bounded in terms of the dimension of $A$: $\pi _0(A)\le 2^{\dim
(A)-1}$. That is, abelian varieties whose ring of endomorphisms is
an order in a totally real number field admit several but not {\em
arbitrarily many} principal polarizations.

It could be expected that Lange's or some other bound for $\pi
_0(A)$ held for any simple abelian variety. Hence the question:
{\em given $g\ge 1$, are there simple abelian varieties of
dimension $g$ with arbitrarily many non-isomorphic principal
polarizations?}

As was already observed, this is not the case in dimension $1$. In
$g=2$, only simple abelian surfaces with at most $\pi _0(A)=2$
were known. One of our main results, stated in a particular case,
is the following.

\begin{theorem}
\label{main}

Let $F$ be a totally real number field of degree $[F:\Q ]=n$, let
$R_F$ denote its ring of integers and $\vartheta _{F/\Q }$ the
different of $F$ over $\Q $. Let $A$ be a complex abelian variety
of dimension $2 n$ whose ring of endomorphisms $\End (A)\simeq
\cO$ is a maximal order in a totally indefinite quaternion
division algebra $B$ over $F$.

Assume that the narrow class number $h_+(F)$ of $F$ is $1$ and
that $\vartheta _{F/\Q }$ and $ \disc (\cO )$ are coprime ideals.
Then,

\begin{enumerate}
\item
$A$ is principally polarizable.

\item
The number of isomorphism classes of principal polarizations on
$A$ is

$$
\pi _0(A)= \frac {1}{2} \sum _S h(S),
$$
where $S$ runs among the set of orders in the CM-field $F(\sqrt
{-D})$ that contain $R_F[\sqrt {-D}]$, the element $D\in F^*_+$ is
taken to be a totally positive generator of the reduced
discriminant ideal $\cD _{\cO }$ of $\cO $ and $h(S)$ denotes its
class number.

\end{enumerate}

In particular, if $A$ is an abelian surface,

$$
\pi _0(A)=
  \begin{cases}
    \dfrac {h(-4 D)+h(-D)}{2} &
\text{ if }D\equiv 3\quad \mathrm{ mod }4, \\
    \dfrac {h(-4 D)}{2} & \text{otherwise.}
  \end{cases}
$$

\end{theorem}

We prove Theorem \ref{main} in the more general form of
Proposition \ref{polpral} and our main Theorem \ref{pio}. In order
to accomplish it, we present an approach to the problem which
stems from Shimura's classical work \cite{Sh1} on analytic
families of abelian varieties with prescribed endomorphism ring.

Our approach is essentially different from Lange's in \cite{La} or
Ibukiyama-Katsura-Oort's in \cite{IbKaOo}. Indeed, whereas in
\cite{La} and \cite{IbKaOo} the (noncanonical) interpretation of
line bundles as symmetric endomorphisms is exploited, we translate
the questions we are concerned with to Eichler's language of
optimal embeddings. This leads us to solve a problem that has its
roots in the work of O'Connor, Pall and Pollack (cf.\,\cite{Po})
and that has its own interest: see Section \ref{poll} for details.

In regard to the question above, the second main result of this
article is the following.

\begin{theorem}
\label{ass}

Let $g$ be a positive integer. Then
\begin{enumerate}
\item
If $g$ is even, there exist simple abelian varieties $A$ of
dimension $g$ such that $\pi _0(A)$ is arbitrarily large.
\item
If $g$ is odd and square-free, $\pi _0(A)\leq 2^{g-1}$ for any
simple abelian variety $A$ of dimension $g$ over $\C $.

\end{enumerate}

\end{theorem}

The boundless growth of $\pi _0(A)$ when $g$ is even follows from
our main Theorem \ref{pio} combined with analytical results on the
asymptotic behaviour and explicit bounds for relative class
numbers of CM-fields due to Horie-Horie (\cite{HoHo}) and
Louboutin (\cite{Lo1}, \cite{Lo2}). The second part of Theorem
\ref{ass} follows from the ideas of Lange in \cite{La}. See
Section \ref{Asymp} for details.

The following corollary follows from Theorem \ref{ass} and the
fact that any simple principally polarized abelian surface is the
Jacobian of a smooth curve of genus $2$ which, by Torelli's
Theorem, is unique up to isomorphism.

\begin{corollary}
\label{jac}

There are arbitrarily large sets $C_1$ , ..., $C_N$ of pairwise
noniso\-mor\-phic genus $2$ curves with isomorphic simple
unpolarized Jacobian varieties $J(C_1)\cong J(C_2)\cong ... \cong
J(C_N)$.

\end{corollary}

In view of Theorem \ref{ass}, it is natural to wonder whether
there exist arbitrarily large sets of pairwise nonisomorphic
curves of given even genus $g\geq 4$ with isomorphic unpolarized
Jacobian varieties. In this direction, Ciliberto and van der Geer
(\cite{CivdGe}) proved the existence of two nonisomorphic curves
of genus 4 whit isomorphic Jacobian varieties. Explicit examples
of curves with isomorphic (nonsimple) Jacobians have been
constructed by Howe (\cite{Ho}), while examples of pairs of
distinct modular curves of genus $2$ defined over $\Q $ with
isomorphic unpolarized absolutely simple Jacobian varieties have
been obtained by Gonz\'{a}lez, Gu\`{a}rdia and the author in
\cite{GoGuRo}.

Finally, let us note that the statement of Theorem \ref{ass} does
not cover odd non square-free dimensions.

\begin{conjecture}

Let $g$ be a non square-free positive integer. Then there exist
simple abelian varieties of dimension $g$ such that $\pi _0(A)$ is
arbitrarily large.

\end{conjecture}

The conjecture is motivated by the fact that, when $g$ is odd and
non square-free, there exist abelian varieties whose ring of
endomorphisms is an order in a non commutative division algebra
over a CM-field and there is a strong similitude between the
arithmetic of the N\'{e}ron-Severi groups of these abelian
varieties and those in the quaternion case.

{\em Acknowledgements. } I am indebted to Pilar Bayer for her
assistance throughout the elaboration of this work. I also express
my gratitude to J. C. Naranjo, S. Louboutin, J. Brzezinski and G.
van der Geer for some helpful conversations. I thank N.
Schappacher and the Institut de Rech\'{e}rche Math\'{e}matique
Avanc\'{e}e at Strasbourg for their warm hospitality in March
2001. Finally, I thank the referee for the valuable help in
improving the exposition.

\section{{\bfseries
Abelian varieties with quaternionic multiplication and their
N\'{e}ron-Severi group}}

Let $F$ be a totally real number field of degree $[F:\Q ]=n$ and
let $R_F$ be its ring of integers. Let $B$ denote a totally
indefinite division quaternion algebra over $F$ and let $\cD =
\disc (B) = \prod _{i=1}^{2 r}\wp _{i}$, where $\wp _i$ are finite
prime ideals of $F$ and $r\geq 1$, be its (reduced) discriminant
ideal. We shall denote $\n =\n _{B/F}$ and $\tr = \tr _{B/F}$ the
reduced norm and trace of $B$ over $F$, respectively. Since $B$ is
totally indefinite, the Hasse-Schilling-Maass Norm Theorem asserts
that $n(B)=F$ (cf.\,\cite{HaSc} and \cite{Vi}, p.\,80). We fix an
isomorphism of $F$-algebras

$$
(\eta _{\sigma }): B\otimes _{\Q }\R \simeq \bigoplus _{\sigma
}\mbox{M}_2(\R ^{\sigma }),
$$
where $\sigma :F\hookrightarrow \R $ runs through the set of
embeddings of $F$ into $\R $ and $\R ^{\sigma }$ denotes $\R $ as
a $F$-vector space via the immersion $\sigma $. For any $\beta \in
B$, we will often abbreviate $\beta ^{\sigma }=\eta _{\sigma
}(\beta )\in \mbox{M}_2(\R )$.

Let $\cO $ be an hereditary order of $B$, that is, an order of $B$
all whose one-sided ideals are projective. The discriminant ideal
$\cD _{\cO }=\disc (\cO )$ of $\cO $ is square-free and it can be
written as $\cD \cdot \cN $ for some ideal $\cN $ coprime to $\cD
$, which is called the {\em level} of $\cO $ (\cite{Re}, Chapter
9).

\begin{definition}

An abelian variety $A$ has quaternionic multiplication by $\cO $
if $\dim (A) = 2 n$ and $\End (A)\simeq\cO $.

\end{definition}

For the rest of this section, let $A/\C $ denote a complex abelian
variety with quaternionic multiplication by the hereditary order
$\cO $. We will identify $\End (A)=\cO $ and $\End (A)\otimes \Q
=B$. Since $B$ is a division algebra, $A$ is simple, that is, it
contains no proper sub-abelian varieties.

As a complex manifold, $A(\C )=V/\Lambda $ for $V$ a complex
vector space of dimension $g$ and $\Lambda \subset V$ a co-compact
lattice that can be identified with the first group of integral
singular homology $H_1(A,\Z )$. The lattice $\Lambda $ is
naturally a left $\cO $-module and $\Lambda \otimes \Q $ is a left
$B$-module of the same rank over $\Q $ as $B$. Since every left
$B$-module is free (cf.\,\cite{We}, Chapter 9), there is an
element $v_0\in V$ such that $\Lambda \otimes \Q = B\cdot v_0$ and
therefore $\Lambda = \cI \cdot v_0$ for some left $\cO $-ideal
$\cI \subset B$.

Let $\Pic _{\ell }(\cO )$ be the pointed set of left (projective)
ideals of $\cO $ up to principal ideals. By a theorem of Eichler
(\cite{Ei1}, \cite{Ei2}), the reduced norm on $B$ induces a
bijection of sets $\n :\Pic _{\ell } (\cO )\stackrel{\sim }{\ra
}\Pic (F)$ onto the class group of $F$.

Note that the left ideal $\cI $ is determined by $A$ up to
principal ideals and we can choose (and fix) a representative of
$\cI $ in its class in $\Pic _{\ell }(\cO )$ such that
$\mathrm{n}(\cI )\subset F$ is coprime with $\cD _{\cO }$. This is
indeed possible because $B$ is totally indefinite: it is a
consequence of the Hasse-Schilling-Maass Norm Theorem, Eichler's
Theorem quoted above and the natural epimorphism of ray class
groups $\mathrm{Cl}^{\cD _{\cO }}(F)\ra \mathrm{Cl}(F)$ of ideals
of $F$ (\cite{Ne}, Chapter VI, Section 6).

Let $\rho _a: B\hookrightarrow \End (V)\simeq \mbox{M}_{2n}(\C )$
and $\rho _r: \cO \hookrightarrow \End (\Lambda )\simeq
\mbox{M}_{4n}(\Z )$ denote the analytic and rational
representations of $B$ and $\cO $ on $V$ and $\Lambda $,
respectively. It is well known that $\rho _r\sim \rho _a \bigoplus
\bar \rho _a$ and it follows that, in an appropriate basis,

$$
\rho _a(\beta )=\mathrm{ diag }(\eta _{\sigma _i}(\beta))
$$
for any $\beta \in B$ (cf.\,\cite{LaBi}, Chapter 9, Lemma 1.1).
Moreover, this basis can be chosen so that the coordinates of
$v_0$ are $(\tau _1, 1, ..., \tau _n, 1)$ for certain $\tau _i\in
\C$, Im$(\tau _i)\not =0$. The choice of the element $v_0$ fixes
an isomorphism of real vector spaces $B\otimes \R \simeq V$.

Conversely, for any choice of a vector $v_0 = (\tau _1, 1, ...,
\tau _n, 1)\in V$ with $\mathrm{Im }(\tau _i) \not = 0$ and a left
$\cO $-ideal $\cI $ in $B$, we can consider the complex torus
$V/\Lambda $ with $\Lambda = \cI \cdot v_0$ and $B$ acting on $V$
via the fixed diagonal analytic representation $\rho _a$. The
torus $V/\Lambda $ admits a polarization and can be embedded in a
projective space. In consequence, it is the set of complex points
of an abelian variety $A$ such that $\End (A)\supseteq \cO $
(cf.\,\cite{Sh1} and \cite{LaBi}, Chapter 9, Section 4).

If $\cO =\cO _{\ell } (\cI )=\{ \beta \in B, \beta \cI \subseteq
\cI  \} $ is the left order of $\cI $ in $B$, it holds that for
the choice of $v_0$ in a dense subset of $V$, we exactly have
$\End (A) = \cO $. Besides, for $v_0$ in a subset of measure zero
of $V$, $A$ fails to be simple and it is isogenous to the square
$A_0^2$ of an abelian variety of dimension $n$ such that $\End
(A_0)$ is an order in a purely imaginary quadratic extension of
$F$ (\cite{Sh2}, Section 9.4).

Let $\NS (A)=\mathrm{Pic}(A)/\mathrm{Pic}^0(A)$ be the
N\'{e}ron-Severi group of line bundles on $A$ up to algebraic
equivalence. Two line bundles $\cL _1$ and $\cL _2\in
\mathrm{NS}(A)$ are said to be isomorphic, denoted $\cL _1\simeq
\cL _2$, if there is an automorphism $\alpha \in \Aut (A)$ such
that $\cL _1=\alpha ^*(\cL _2)$.

For an $\cO $-left ideal $J$, let

$$
J ^{\sharp }=\{\beta \in B: \tr _{B/\Q }(J \beta )\subseteq \Z \}
$$
be the codifferent of $J $ over $\Z $. It is a right ideal of $\cO
$ projective over $R_F$. If we let $B_0=\{ \beta \in B: \tr (\beta
)=0\}$ denote the additive subgroup of pure quaternions of $B$, we
put $J _0^{\sharp }=J ^{\sharp } \cap B_0$. Finally, we define
$\cN (J )=\n (J ) \cO =J \bar J $ to be the two-sided ideal of
$\cO $ generated by the ideal $\n (J )$ of $F$ of reduced norms of
elements in $J $.

The following theorem describes $\mathrm{NS}(A)$ intrinsically in
terms of the arithmetic of $B$. In addition, it establishes when
two line bundles on $A$ are isomorphic and translates this into a
certain conjugation relation in $B$. We keep the notations as
above.

\begin{theorem}
\label{NeSe} There is a natural isomorphism

$$
\begin{array}{clc}
\mathrm{NS}(A)&\stackrel {\sim }{\rightarrow } & \cN (\cI )_0^{\sharp } \\
\cL &\mapsto &\mu =\mu (\cL )
\end{array}
$$
between the N\'{e}ron-Severi group of $A$ and the group of pure
quaternions of the codifferent of the two-sided ideal $\cN (\cI
)$. Moreover, for any two line bundles $\cL _1$ and $\cL _2\in
\mathrm{NS}(A)$, we have that $\cL _1\simeq \cL _2$ if and only if
there exists $\alpha \in \cO ^*$ such that $\mu (\cL _1) = \bar
\alpha \mu (\cL _2) \alpha $.

\end{theorem}

{\em Proof. } By the Appell-Humbert Theorem, the first Chern class
allows us to interpret a line bundle $\cL \in \mathrm{NS}(A)$ as a
Riemann form: an $\R $-alternate bilinear form $E=c_1(\cL
):V\times V\rightarrow \R $ such that $E(\Lambda \times \Lambda
)\subset \Z $ and $E(\sqrt {-1} u,\sqrt {-1} v)=E(u,v) $ for all
$u$, $v\in V$. Fix a line bundle $\cL $ on $A$ and let $E=c_1(\cL
)$ be the corresponding Riemann form. The linear map $B\rightarrow
\Q $, $\beta \mapsto E(\beta v_0, v_0)$ is a trace form on $B$ and
hence, by the nondegeneracy of $\tr _{B/\Q }$, there is a unique
element $\mu \in B$ such that $E(\beta v_0, v_0)=
\mathrm{tr}_{B/\Q }(\mu \beta )$ for any $\beta \in B$. Since $E$
is alternate, $E(a v_0,a v_0)= \mathrm{tr}_{F/\Q }(a^2
\mathrm{tr}_{B/F}(\mu ))=0$ for any $a\in F$. It follows again
from the nondegeneracy of $\mathrm{tr}_{F/\Q }$ and the fact that
the squares $F^{*2}$ span $F$ as a $\Q $-vector space that
$\mathrm{tr}_{B/F}(\mu )=0$. Thus $\mu ^2+\delta =0$ for some
$\delta \in F$.

The line bundle $\cL $ induces an anti-involution $\varrho $ on
$B$ called the Rosati involution. It is characterized by the rule
$E(u, \beta v)=E(\beta ^{\varrho }u, v)$ for any $\beta \in B$ and
$u$, $v\in V$. From the discussion above, it must be $\beta
^{\varrho }=\mu ^{-1}\bar \beta \mu $ and we conclude that the
Riemann form $E=c_1(\cL )$ attached to the line bundle $\cL $ on
$A$ is

$$
\begin{array}{cccc}
E:=E_{\mu }: & V\times V & \longrightarrow & \R \\
  & (u,v) & \mapsto &\tr _{B\otimes _{\Q }\R /\R }(\mu \gamma \bar \beta )
\end{array}
$$
where $\mu \in B$ is determined as above and $\gamma $ and $\beta
$ are elements in $B\otimes _{\Q }\R \simeq \M _2(\R )^n$ such
that $u=\gamma v_0$ and $v=\beta v_0$. Since $E(\Lambda \times
\Lambda ) \subset \Z $ and $\tr (\mu )=0$, we deduce that $\mu \in
\cN (\cI )_0^{\sharp }$.

Conversely, one checks that any element $\mu \in \cN (\cI
)_0^{\sharp }$ defines a Riemann form $E_{\mu }$ which is in turn
the first Chern class of a line bundle $\cL $ on $A$. Indeed,
since $\mu \in \cN (\cI )^{\sharp }$, the form $E_{\mu }$ is
integral over the lattice $\Lambda =\cI \cdot v_0$ and $E_{\mu }$
is alternate because $\mathrm{tr}(\mu )=0$. Moreover, let $\iota
=\mathrm{diag}(\iota _1, ..., \iota _n)\in \mathrm{GL} _{2n}(\R )$
with $\iota _i\in \mathrm{GL} _2(\R )$ and $\iota _i^2+1=0$, be a
matrix such that $\iota \cdot v_0=\sqrt{-1} v_0$. Then $E_{\mu
}(\sqrt{-1}u, \sqrt{-1}v) = E_{\mu }(\gamma \sqrt{-1}v_0, \beta
\sqrt{-1}v_0) = \mathrm{tr} (\mu \gamma \iota \bar \iota \bar
\beta ) = E_{\mu }(u, v)$ for all $u$, $v\in V$. This concludes
the proof of the first part of the theorem.

As for the second, we note that the first Chern class of the
pull-back $\alpha ^*\cL $ of a line bundle $\cL $ on $A$ by an
automorphism $\alpha \in \Aut (A)=\cO ^*$ is represented by the
Riemann form $\alpha ^*E: V\times V\rightarrow \R $ defined by
$(u, v)\mapsto E(\alpha u, \alpha v)$, where $E=c_1(\cL )$ is the
Riemann form associated to $\cL $. Hence, if $\cL _1=\alpha ^*
(\cL _2)$, then tr$(\mu _1 \gamma \bar \beta )=\mathrm{tr}(\mu _2
\alpha \gamma \bar \beta \bar \alpha )=\mathrm{tr}(\bar \alpha \mu
_2 \alpha \gamma \bar \beta )$ for all $\gamma $ and $\beta \in B$
and this is satisfied if and only if $\mu _1=\bar \alpha \mu
_2\alpha $, by the nondegeneracy of the trace form. Reciprocally,
one checks that if $\mu _1=\bar \alpha \mu _2\alpha $ for some
$\alpha \in \cO ^*$, then $E_{\mu _1}=\alpha ^*E_{\mu _2}$ and
therefore $\cL _1=\alpha ^*\cL _2$. $\Box $

In view of Theorem \ref{NeSe}, we identify the first Chern class
$c_1(\cL )$ of a line bundle $\cL $ on $A$ with the quaternion
$\mu =\mu (\cL ) \in B_0$ such that $E_{\mu }$ is the Riemann form
associated to $\cL $. We are led to introduce the following
equivalence relation, that was studied (over $B$) by O'Connor and
Pall in the 1930s and Pollack in the 1960s (cf.\,\cite{Po}).

\begin{definition}\label{pollconj} Two
quaternions $\mu _1$ and $\mu _2\in B$ are {\em Pollack
conjugated} over $\cO $ if $\mu _1 = \bar \alpha \mu _2 \alpha $
for some unit $\alpha \in \cO ^*$. We will denote it by $\mu
_1\sim _p \mu _2$.
\end{definition}

\section{{\bfseries
Isomorphism classes of line bundles and Eichler theory on optimal
embeddings}}\label{proofpi}

As in the previous section, let $A$ denote an abelian variety with
quaternionic multiplication by an hereditary order $\cO $ in a
totally indefinite division quaternion algebra $B$ over a totally
real number field $F$. As is well-known, a line bundle $\cL \in
\mathrm{NS}(A)$ induces a morphism $\varphi _{\cL }:A\rightarrow
\hat A $ defined by $P\mapsto \mbox{t}_P^*(\cL )\otimes \cL ^{-1}
$, where $t_P:A\rightarrow A$ denotes the translation by $P$ map.
Since $A$ is simple, any nontrivial line bundle $L\in
\mathrm{NS}(A)$ is nondegenerate: $\varphi _{\cL }$ is an isogeny
with finite kernel $K(\cL )$. We say that $\cL $ is {\em
principal} if $K(\cL )$ is trivial, that is, if $\varphi _{\cL
}:A\rightarrow \hat A$ is an isomorphism.

\begin{proposition}
\label{degree}

Let $\cL $ be a line bundle on $A$ and let $c_1(\cL ) = \mu $ be
its first Chern class for some element $\mu \in B$ such that $\mu
^2+\delta =0$ and $\delta \in F$. Then

$$
\mathrm{deg }(\varphi _{\cL })={\mathrm{N}_{F/\Q }( \vartheta
_{F/\Q }^2\cdot \n (\cI )^2\cdot \cD _{\cO }\cdot \delta )}^2
$$
where $\vartheta _{F/\Q }=(R_F^{\sharp })^{-1}$ is the different
of $F$ over $\Q $.

\end{proposition}

{\em Proof.} The degree $|K(\cL )|$ of $\varphi _{\cL }$ can be
computed in terms of the Riemann form as follows: $\mathrm{deg
}(\varphi _{\cL }) = \det (E_{\mu }(x_i, x_j)) = \det
(\mathrm{tr}_{B/\Q } (\mu \beta _i \bar \beta _j))$, where
$x_i=\beta _i v_0$ runs through a $\Z $-basis of the lattice
$\Lambda $. We have $\mathrm{det }(\mathrm{ tr}_{B/\Q }(\mu \beta
_i \bar \beta _j))=\mathrm{n}_{B/\Q }(\mu )^2\det (\mathrm{
tr}_{B/\Q } (\beta _i \cdot \bar \beta _j))= \mathrm{n}_{B/\Q
}(\mu )^2 \disc _{B/\Q }(\cI )^2 = (\mathrm{N}_{F/\Q }(\delta
\cdot \mathrm{n}(\cI )^2 \cdot \vartheta _{F/\Q }^2\cdot \cD _{\cO
}))^2. \quad \Box $

As a consequence, we obtain the following criterion that
establishes whether the abelian variety $A$ admits a principal
line bundle in terms of the arithmetic of the hereditary order
$\cO =\End (A)$ in $B$ and the left ideal $\cI \cong H_1(A,\Z )$.
Crucial in the proof of the theorem below is the classical theory
of Eichler optimal embeddings.

For the sake of simplicity and unless otherwise stated, we assume
for the rest of the article that $\vartheta _{F/\Q } \mbox{ and }
\cD _{\cO }\mbox{ are coprime ideals in }F$. The general case can
be dealt by means of the remark below.

\begin{theorem}
\label{pral}

The abelian variety $A$ admits a {\em principal} line bundle if
and only if the ideals $\cD _{\cO }$ and $\vartheta _{F/\Q }\cdot
\mathrm{n}(\cI )$ of $F$ are {\em principal}.

\end{theorem}

{\em Proof. } Let $\cL $ be a principal line bundle on $A$ and let
$E_{\mu }=c_1(\cL )$ be the associated Riemann form for some $\mu
\in \cN (\cI )_0^{\sharp }$ such that $\mu ^2+\delta =0$. Since
$\cL $ is principal, the induced Rosati involution $\varrho $ on
$\End (A)\otimes \Q =B$ must also restrict to $\End (A)=\cO$ and
we already observed that ${\beta }^{\varrho } = \mu ^{-1}\bar
\beta \mu $. Therefore $\mu $ belongs to the normaliser group
$\mathrm{Norm}_{B^*}(\cO )$ of $\cO $ in $B$. The quotient
$\mathrm{Norm}_{B^*}(\cO )/\cO ^*F^* \cong W$ is a finite abelian
2-torsion group and representatives $w$ of $W$ in $\cO $ can be
chosen so that the reduced norms $\n (w)\in R_F$ are only
divisible by the prime ideals $\wp \vert \cD _{\cO }$
(cf.\,\cite{Vi}, p.\,39, 99, \cite{Br}). Hence, we can express
$\mu =u\cdot t \cdot w^{-1}$ for some $u\in \cO ^*$, $t \in F^*$
and $w\in W$.

Recall that $(\n (\cI ), \cD _{\cO })=1$ and $(\vartheta _{F/\Q },
\cD _{\cO })=1$. Since, from Proposition \ref{degree},
$\mathrm{n}(\cI )^{2} \cdot \vartheta _{F/\Q }^{2}\cdot \cD _{\cO
} = (\delta ^{-1}) = (t ^{-2}\cdot \mathrm{n}(w))$, we conclude
that $\mathrm{n}(\cI )\cdot \vartheta _{F/\Q }=(t ^{-1})$ and $\cD
_{\cO }=(\mathrm{n}(w))$ are principal ideals.

Conversely, suppose that $\mathrm{n}(\cI )\cdot \vartheta _{F/\Q
}=(t ^{-1})$ and $\cD _{\cO }=(D)$ are principal ideals, generated
by some elements $t $ and $D\in F^*$. Let $S$ be the ring of
integers in $L = F(\sqrt {-D})$. Since any prime ideal $\wp \vert
D$ ramifies in $L$, Eichler's theory of optimal embeddings
guarantees the existence of an embedding $\iota :S\hookrightarrow
\cO $ of $S$ into the quaternion order $\cO $ (cf.\,\cite{Vi},
p.\,45). Let $w=\iota (\sqrt {-D})\in \cO $ and let $\mu = t \cdot
w^{-1}$. As one checks locally, $\mu\in \mathrm{Norm}_{B^*}(\cO
)\cap \cN (\cI )_0^{\sharp }$ and, by Theorem \ref{NeSe} and
Proposition \ref{degree}, $\mu $ is the first Chern class of a
principal line bundle on A. $\Box $

\begin{corollary}
\label{dual}

If $\cD _{\cO }$ and $\vartheta _{F/\Q }\cdot \mathrm{n}(\cI )$
are principal ideals, then $A$ is self-dual, that is, $A\cong \hat
A$.

\end{corollary}

\begin{remark} The case when $\vartheta _{F/\Q }$ and $\cD _{\cO }$ are non
necessarily coprime is reformulated as follows: $A$ admits a
principal line bundle if and only if there is an integral ideal
$\mathfrak a=\wp _1^{e_1}\cdot ...\cdot \wp _{2r}^{e_{2r}}\vert
\vartheta _{F/\Q }$ in $F$ such that both $\cD _{\cO } \cdot
\mathfrak a^2$ and $\mathrm{n}(\cI )\cdot \vartheta _{F/\Q }\cdot
\mathfrak a^{-1}$ are principal ideals. In this case, $A$ is also
self-dual. The proof is {\em mutatis mutandi} the one given above.
\end{remark}

\begin{definition}

The set of isomorphism classes of principal line bundles on $A$ is
$$
\Pi (A)=\{ \cL \in \NS (A): \deg (\cL )=1\} /_{\simeq }.
$$

\end{definition}

\begin{definition}

Assume that $\cD _{\cO }$ is a principal ideal of $F$. Then, we
let

$$
\cP (\cO )=\{ \mu \in \cO : \mathrm{tr}(\mu )=0, \n (\mu ) R_F =
\cD _{\cO }\}
$$
and we define $P(\cO ) = \cP (\cO )/_{\sim _p}$ to be the
corresponding set of Pollack conjugation classes.

\end{definition}

The above proof, together with Theorem \ref{NeSe}, yields

\begin{corollary}
\label{bij}

Let $A$ be an abelian variety with quaternionic multiplication by
a maximal order $\cO $. If $\cD _{\cO }=(D)$ and $\vartheta _{F/\Q
} \cdot \mathrm{n}(\cI )=(t ^{-1})$ are principal ideals, the
assignation $$\cL \mapsto t \cdot c_1(\cL )^{-1}$$ induces a
bijection of sets between $\Pi (A)$ and $P (\cO )$.

\end{corollary}

In view of Corollary \ref{bij}, it is our aim to compute the
cardinality $\pi (A)= |\Pi (A)|$ of the set of isomorphism classes
of principal line bundles on an abelian variety $A$ with
quaternionic multiplication by a maximal order $\cO $. Theorem
\ref{pi} below exhibits a close relation between $\pi (A)$ and the
class number of $F$ and of certain orders in quadratic extensions
$L/F$ that embed in $B$.

Assume that there is a principal line bundle on $A$. Otherwise
$\pi (A)=0$ and there is nothing to compute. We assume also that
$(\vartheta _{F/\Q }, \cD _{\cO }) = 1$. By Theorem \ref{pral}, we
know that $\cD _{\cO }= (D)$ for some $D\in F^*$. We have

\begin{theorem}
\label{pi}

Let $A$ be an abelian variety with quaternionic multiplication by
a maximal order $\cO $. Then

$$
\pi (A) = \frac {1}{2 h(F)}\sum _u \sum _{S} 2^{e_S} h(S),
$$
where $u\in R_F^*/R_F^{*2}$ runs through a set of representatives
of units in $R_F$ up to squares and $S$ runs through the (finite)
set of orders in $F(\sqrt {-u D})$ such that $R_F[ \sqrt {-u D}]
\subseteq S$. Here, $2^{e_S} = |R_F^*/N_{F(\sqrt {-u
D})/F}(S^*)|$.

\end{theorem}

The proof of Theorem \ref{pi} will be completed in Section
\ref{poll}. There are several remarks to be made for the sake of
its practical applications.

\begin{remark} Note that $R_F^{*2}\subseteq N_{L/F}(S^*)$ and hence
$R_F^*/N_{L/F}(S^*)$ is naturally an $\F _2$-vector space. By
Dirichlet's Unit Theorem, $e_S\leq [F:\Q ] = n$. The case $F=\Q $
is trivial since $\Z ^*/\Z ^{*2} = \{ \pm 1\} $. In the case of
real quadratic fields $F$, explicit fundamental units $u\in R_F^*$
such that $R_F^*/R_F^{*2} = \{ \pm 1, \pm u\} $ are well known.
For totally real number fields of higher degree there is abundant
literature on systems of units. See \cite{La} and \cite{Wa},
Chapter 8 for an account.
\end{remark}

\begin{remark}\label{rem} Let $2 R_F = \mathfrak q_1^{e_1}\cdot ... \cdot \mathfrak
q_m^{e_m}$ be the decomposition of $2$ into prime ideals in $F$.
Fix a unit $u\in R_F^*$. Then, the conductor $\mathfrak f$ of
$R_F[\sqrt {-u D}]$ in $L$ over $R_F$ is $\mathfrak f = \prod
_{\stackrel {\mathfrak q\vert 2}{\mathfrak q\nmid D}} \mathfrak
q^{a_{\mathfrak q}}, \mbox{ for some }0\leq a_{\mathfrak q}\leq
e_{\mathfrak q}$.

Further, the conductor $\mathfrak f$ can be completely determined
in many cases as follows. For a prime ideal $\mathfrak q\vert 2$
such that $\mathfrak q\nmid D$, let $\pi $ be a local uniformizer
of the completion of $F$ at $\mathfrak q$ and $k = \F _{2^f}$ be
the residue field. Let $e = e_{\mathfrak q}\geq 1$. Since $-u D\in
R_{F_{\mathfrak q}}^*$, we have $-u D = x_0 + x_k \pi ^{k} +
x_{k+1} \pi ^{k+1} + ... $ for some $1\leq k\leq \infty$ and $x_i$
in a system of representatives of $\F _{2^f}$ in $R_{F_{\mathfrak
q}}$ such that $\bar {x_0}$ and $\bar {x_k}\not = 0$. Here we
agree to set $k=\infty $ if $-u D = x_0$. Then, $\mathrm{min
}([\frac {k}{2}], e)\leq a_{\mathfrak q}\leq e$ and we exactly
have

$$
a_{\mathfrak q} = \begin{cases}
    [\frac {k}{2}] & \text{if } k\leq e+1,\\
    e & \text{if }[\frac {k}{2}]\geq e.
  \end{cases}
$$

Otherwise, if $[\frac{k}{2}]<e<k-1$, the determination of
$a_{\mathfrak q}$ depends on the choice of the system of
representatives of $\F _{2^f}$ in $R_{F_{\mathfrak q}}$ and it
deserves a closer inspection.

This gives an easy criterion for deciding whether $R_F[\sqrt {-u
D}]$ is the ring of integers of $F(\sqrt {-u D})$ (that is,
$\mathfrak f = 1$). In any case, the set of orders $S$ in $L =
F(\sqrt {-u D})$ that contain $R_F[\sqrt {-u D}]$ can be described
as follows. Any order $S\supseteq R_F[\sqrt {-u D}]$ has conductor
$\mathfrak f _S\vert \mathfrak f $ and for every ideal $\mathfrak
f '\vert \mathfrak f $ there is a unique order $S\supseteq
R_F[\sqrt {-u D}]$ of conductor $\mathfrak f '$. Further,
$\mathfrak f _S\vert \mathfrak f _T$ if and only if $S\supseteq
T$. We omit the details of the proof of these facts.
\end{remark}

In order to prove Theorem \ref{pi}, we begin by an equivalent
formulation of it. As it was pointed out in Corollary \ref{bij},
the first Chern class induces a bijection of sets between $\Pi
(A)$ and the set of Pollack conjugation classes $P (\cO )$. For
$u\in R_F^*$, let us write $\cP (u, \cO ) :=\{ \mu \in \cO : \mu
^2 + u D = 0\} $. Observe that $\cP (\cO )$ is the disjoint union
of the sets $\cP (u_k, \cO )$ as $u_k$ run among units in any set
of representatives of $R_F^*/R_F^{*2}$.

Any quaternion $\mu \in \cP (u, \cO )$ induces an embedding

$$
\begin{array}{clc}
i_{\mu }: F(\sqrt {-u D})&\hookrightarrow &B \\ a+b\sqrt{-u
D}&\mapsto &a+b\mu
\end{array}
$$
for which $i_{\mu }(R_F[\sqrt {-u D}])\subset \cO $. The following
definition is due to Eichler.

\begin{definition}
Let $S$ be an order over $R_F$ in a quadratic algebra $L$ over
$F$. An embedding $i:S\hookrightarrow \cO $ is {\em optimal} if
$i(S)=i(L)\cap \cO $.
\end{definition}

For any $\mu \in \cP (u, \cO )$ there is a uniquely determined
order $S_{\mu }\supseteq R_F[\sqrt {-u D}]$ such that $i_{\mu }$
is optimal at $S_{\mu }$. Moreover, two equivalent quaternions
$\mu _1\sim _p\mu _2\in \cP (u, \cO )$ are optimal at the same
order $S$. Indeed, if $\alpha \in \cO^*$ is such that $\mu _1=\bar
\alpha \mu _2 \alpha $, then $\alpha $ is forced to have reduced
norm $\mathrm{n}(\alpha )=\pm 1$. Hence $\bar \alpha =\pm \alpha
^{-1} \in \cO ^*$ and the observation follows since $\alpha $
normalizes $\cO $. Conversely, any optimal embedding
$i:S\hookrightarrow \cO $, $S\supseteq R_F[\sqrt {-u D}]$
determines a quaternion $\mu =i(\sqrt {-u D})\in \cP (u, \cO )$.

For any quadratic order $S$ over $R_F$, let $\cP (S, \cO )$ denote
the set of optimal embeddings of $S$ in $\cO $ and $P(S, \cO )=
\cP (S, \cO )/_{\sim _p}$. We obtain a natural identification of
sets

$$
P (\cO ) = \sqcup _{k} \sqcup _{S} P (S, \cO ),
$$
where $S$ runs through the set of quadratic orders $S\supseteq
R_F[\sqrt {-u_k D}]$ for any unit $u_k$ in a set of
representatives of $R_F^*/R_F^{*2}$. Hence, in order to prove
Theorem \ref{pi}, it is enough to show that, for any quadratic
order $S \supseteq R_F[\sqrt {-u D}]$, $u\in R_F^*$, it holds that

$$
p(S, \cO ):=|P (S, \cO )|=\frac {2^{e_S-1} h(S)}{h(F)}.
$$
Since this question is interesting on its own, we will make our
statements in greater generality in the next section.

\section{{\bfseries Pollack conjugation versus Eichler conjugation}}
\label{poll}

Let $F$ be a number field and let $B$ be a division quaternion
algebra over $F$. In \cite{Po}, Pollack studied the obstruction
for two pure quaternions $\mu _1$ and $\mu _2\in B$ with the same
reduced norm to be conjugated over $B^*$, that is, $\mu _1=\bar
{\alpha }\mu _2\alpha $ with $\alpha \in B^*$, and expressed it in
terms of the $2$-torsion subgroup $\mathrm{Br}_2(F)$ of the Brauer
group $\mathrm{Br}(F)$ of $F$. He further investigated the
solvability of the equation $\mu _1=\bar {\alpha }\mu _2\alpha $
over $\cO ^*$ for quaternions $\mu _1$ and $\mu _2$ in a maximal
order $\cO $ of $B$.

As a refinement of his considerations, it is natural to consider
the set of orbits of pure quaternions $\mu \in \cO $ of fixed
reduced norm $\n (\mu )=d\in F^*$ under the action of the group of
units $\cO ^*$ by Pollack conjugation. As already mentioned, a
necessary condition for $\mu _1\sim _p\mu _2$ over $\cO ^*$ is
that $\mu _1$ and $\mu _2$ induce an optimal embedding at the same
quadratic order $F(\sqrt {-d})\supset S\supseteq R_F[\sqrt {-d}]$.
We will drop the restriction on $\cO $ to be maximal in our
statements.

The connection between optimal embeddings orbits and class numbers
is made possible by the theory of Eichler. However, in contrast to
Pollack conjugation, two optimal embeddings $i, j:S\hookrightarrow
\cO $ lie on the same conjugation class in the sense of Eichler,
written $i\sim _e j$, if there exists $\alpha \in \cO ^*$ such
that $i = \alpha ^{-1} j \alpha $. We let $E(S, \cO )=\cP (S, \cO
)/_{\sim _e}$ be the set of Eichler conjugation classes of optimal
embeddings of $S$ into an order $\cO $ of $B$ and denote $e(S, \cO
) = |E(S, \cO )|$.

\begin{proposition}\label{SS}

Let $S$ be an order in a quadratic algebra $L$ over $F$ and let
$\cO $ be an order in a division quaternion algebra $B$ over $F$.
Then, the number of Pollack conjugation classes is

$$
p(S, \cO ) = |\n (\cO ^*)/N_{L/F}(S^*)|\cdot \frac {e(S, \cO
)}{2}.
$$

\end{proposition}

{\em Proof. } Let us agree to say that two pure quaternions $\mu
_1$ and $\mu _2\in \cO $ lie in the same $\pm $Eichler conjugation
class if there exists $\alpha \in \cO ^* $ such that $\mu _1 = \pm
\alpha ^{-1}\mu _2 \alpha $. We shall denote it by $\mu _1 \sim
_{\pm e} \mu _2 $ and $E _{\pm } (S, \cO ) = \cP (S, \cO )/_{\sim
_{\pm e}}$. The identity map $\mu \mapsto \mu$ descends to a
natural surjective map
$$
\begin{array}{clc}
\rho :P (S, \cO )&\ra &E _{\pm } (S, \cO )\\
\end{array}
$$
and the proposition follows from the following lemma.

\begin{lemma}
\label{epsilon2}

Let $e_S = \mathrm{dim }_{\F _2}(\n (\cO ^*)/N_{L/F}(S^*))$. Let
$\mu \in \cP (S, \cO )$ and let $\varepsilon _{\mu }=1$ if $\mu
\sim _e -\mu $ and $\varepsilon _{\mu }=2$ otherwise. Then, in the
$\pm $Eichler conjugation class $\{ \pm \alpha ^{-1}\mu \alpha :
\alpha \in \cO ^*\}$ of $\mu $, there are exactly $\varepsilon
_{\mu } 2^{e_S-1}$ Pollack conjugation classes of pure
quaternions.

\end{lemma}

{\em Proof of Lemma \ref{epsilon2}. } Suppose first that
$\varepsilon _{\mu }=1$. Then, the $\pm $Eichler conjugation class
of $\mu \in \cP (S, \cO )$ is $\{\alpha ^{-1}\mu \alpha :\alpha
\in \cO ^*\} $. Let $\gamma \in \cO ^*$ be such that $-\mu =
\gamma \mu \gamma ^{-1} = \gamma ^{-1}\mu \gamma $. We claim that,
for any given $\alpha \in \cO ^*$, it holds that $\mu \sim _p
\alpha ^{-1}\mu \alpha $ if and only if $\n (\alpha )\in
N_{L/F}(S^*)\cup (-\n (\gamma )N_{L/F}(S^*))$. Indeed, if $\mu
\sim _p \alpha ^{-1}\mu \alpha $, let $\beta \in \cO ^*$ with $\n
(\beta )=\pm 1$ be such that $\bar \beta \alpha ^{-1}\mu \alpha
\beta = \mu $. If $\n (\beta)=1$, then $\alpha \beta \mu = \mu
\alpha \beta $ and hence $\alpha \beta \in L\cap \cO ^*=S^*$. Thus
$\n (\alpha \beta )=\n (\alpha )\in N_{L/F}(S^*)$. If $\n (\beta
)=-1$, a similar argument shows that $\n (\alpha )\in -\n (\gamma
)\cdot N_{L/F}(S^*)$. Conversely, let $\n (\alpha ) = v\in
N_{L/F}(S^*)\cup (-\n (\gamma )N_{L/F}(S^*))$ and let $s\in S^*$
be such that $N_{L/F}(s) = v \mbox{ or} -v \n (\gamma )^{-1}$.
Since $\mu $ induces an embedding $S\hookrightarrow \cO $, we can
regard $s$ as an element in $\cO ^*$ such that $\n (s) = v \mbox{
or} -v \n (\gamma )^{-1}$ and $s \mu = \mu s$. Hence $\alpha ^{-1}
\mu \alpha = \alpha ^{-1} s \mu s^{-1} \alpha = (\alpha ^{-1}
s)\mu \overline {(\alpha ^{-1} s)} \mbox{ or } (\alpha ^{-1} s
\gamma ^{-1})\mu \overline {(\alpha ^{-1} s \gamma ^{-1})}$. This
proves the claim.

Since $B$ is division, Pollack's Theorem on Pall's Conjecture
(\cite{Po}, Theorem 4) applies to show that $-\n (\gamma ) \not
\in N_{L/F}(S^*)$. We then conclude that the distinct Pollack
conjugation orbits in $\{\alpha ^{-1}\mu \alpha :\alpha \in \cO
^*\} $ are exactly the classes

$$
\mathcal C _u = \{ \alpha ^{-1}\mu \alpha : \n (\alpha )\in u
N(S^*)\cup (-\n (\gamma ) u N(S^*))\}
$$
as $u\in \n (\cO ^*)$ runs through a set of representatives in $\n
(\cO ^*)/ \langle -\n (\gamma ), N(S^*)\rangle $. There are
$2^{e_S-1}$ of them.

Assume that $\varepsilon _{\mu }=2$, that is, $\mu \not \sim _{e}
-\mu $. Then, the $\pm $Eichler conjugation class of $\mu \in \cP
(S, \cO )$ is $\{\alpha ^{-1}\mu \alpha :\alpha \in \cO ^*\} \cup
\{-\alpha ^{-1}\mu \alpha :\alpha \in \cO ^*\} $. As in the
previous case, it is shown that $\mu \sim _p \alpha ^{-1}\mu
\alpha $ if and only if $\n (\alpha )\in N_{L/F}(S^*)$ and $\mu
\sim _p -\alpha ^{-1}\mu \alpha $ if and only if $\n (\alpha )\in
-N_{L/F}(S^*)$. We obtain that, as $u\in \n (\cO ^*)$ runs through
a set of representatives in $\n (\cO ^*)/N(S^*)$, the $2^{e_S}$
distinct Pollack conjugation classes in the $\pm $Eichler
conjugation class of the quaternion $\mu \in \cP (S, \cO )$ are

$$
\mathcal C '_u = \{ \alpha ^{-1}\mu \alpha :\n (\alpha )\in u
N(S^*)\} \cup \{ -\alpha ^{-1}\mu \alpha :\n (\alpha )\in -u
N(S^*)\}. \quad \Box
$$

{\em Proof of Theorem \ref{pi}. } Firstly, under the assumptions
of Theorem \ref{pi}, the Hasse-Schilling-Maass Norm Theorem in its
integral version (\cite{Vi}, p.\,90) asserts that $\n (\cO
^*)=R_F^*$. Secondly, we have that

$$
e (S, \cO )= \frac {h(S)}{h(F)}
$$
for any $S\supseteq R_F[\sqrt {-u D}]$, $u\in R_F^*$. This follows
from a theorem of Eichler (cf.\,\cite{Vi}, p.\,98) together with
Remark \ref{rem}. The combination of these facts together with the
discussion at the end of the Section \ref{proofpi} and Proposition
\ref{SS} yield the theorem. $\quad \Box $

\begin{remark} In view of Proposition \ref{SS}, the effective computation of the
number of Pollack conjugation classes $p(S, \cO )$ for arbitrary
orders lies on the computability of the groups $N_{L/F}(S^*)$ and
$\n (\cO ^*)$ and the number $e(S, \cO )$. The study of the former
depends on the knowledge of the group of units $S^*$ and there is
abundant literature on the subject. If $\cO $ is an Eichler order,
the Hasse-Schilling-Maass Norm Theorem in its integral version
(\cite{Vi}, p.\,90) describes $\n (\cO ^*)$ in terms of the
archimedean ramified places of $B$. Finally, there are several
manuscripts which deal with the computation of the numbers $e(S,
\cO )$. See \cite{Vi}, p.\,45, for Eichler orders and
\cite{HiPiSh} and \cite{Br} for Gorenstein and Bass orders.
\end{remark}

Since it will be of use later, we will consider in Proposition
\ref{Galois} below a stronger form of Eichler's Theorem. We keep
the notations as in Theorem \ref{pi}. Let $S\supseteq R_F[\sqrt
{-u D}]$, $u\in R_F^*$, be a quadratic order and let $H_S$ be the
{\em ring class field} of $S$ over $L=F(\sqrt {-u D})$. The Galois
group $\mathrm{Gal} (H_S/L)$ is isomorphic, via the Artin
reciprocity map, to the Picard group $\Pic (S)$ of classes of
locally invertible ideals of $S$. In the particular case that $S$
is the ring of integers of $L$, then $H_S$ is the {\em Hilbert
class field} of $L$. The quadratic extension $L/F$ is ramified at
the prime ideals $\wp \vert D$ and recall that, by Remark
\ref{rem}, these prime ideals do not divide the conductor of $S$.
Therefore $L$ and $H_S$ are linearly disjoint over $F$, that is,
$F = L\cap H_S$. The norm induces a map $N_{L/F}:\Pic
(S)\rightarrow \Pic (R_F)$ that, by the reciprocity isomorphism
can be interpreted as the restriction map $\mathrm{Gal} (H_S/L)
\rightarrow \mathrm{Gal} (L\cdot H_F/L)\simeq \mathrm{Gal}
(H_F/F)$ (cf.\,\cite{Ne}, Chapter VI, Section 5). In particular,
we have an exact sequence

$$
0\rightarrow \Delta \rightarrow \Pic (S)\stackrel
{N_{L/F}}{\rightarrow } \Pic (R_F)\rightarrow 0 .
$$
Here, $\Delta = \mbox{Ker }(N_{L/F})$ can be viewed as the Galois
group of $H_S$ over the fixed field $L_{\bigtriangleup }$ of $H_S$
by $\Delta $. The group $\Delta = \mathrm{Gal}
(H_S/L_{\bigtriangleup })$ acts on $E (S, \cO )$ by a reciprocity
law as follows: let $i:S\hookrightarrow \cO $ be an optimal
embedding and let $\tau \in \mathrm{Gal} (H_S/L_{\bigtriangleup
})$. Then let $\mathfrak b= [\tau , H_S/L]$ be the locally
invertible ideal in $S$ corresponding to $\tau $ by the Artin's
reciprocity map. Since the reduced norm on $B$ induces a bijection
of sets $\n :\Pic _{\ell }(\cO )\simeq \Pic (F)$ and
$N_{L/K}(\mathfrak b)$ is a principal ideal in $F$, it follows
that $i(\mathfrak b)\cO =\beta \cO $ is a principal right ideal of
$\cO $ and we can choose a generator $\beta \in \cO $. Then $\tau
$ acts on $i\in E (S, \cO )$ as

$$
i^{\tau } = \beta ^{-1}i \beta .
$$
It can be checked that this action does not depend on the choice
of the ideal $\mathfrak b$ in its class in $\Pic (S)$ nor on the
choice of the element $\beta \in \cO $. Moreover, a local argument
shows that this action is free. Since $|\Delta | = |E (S, \cO )|$,
we obtain

\begin{proposition}
\label{Galois}

The action of $\Delta $ on the set of Eichler conjugacy classes of
optimal embeddings $E (S, \cO )$ is free and transitive.

\end{proposition}

The above action acquires a real arithmetic meaning and coincides
with Shimu\-ra's reciprocity law in the particular case that $L$
is a CM field over $F$. In this situation, $E (S, \cO )$ can also
be interpreted as the set of {\em Heegner points} on a Shimura
variety $\mathfrak X$ on which the Galois group $\Delta $ is
acting (cf.\,\cite{Sh2}, Section 9.10).

\section{{\bfseries The index of a nondegenerate line bundle}}

Let $\cL \in \mathrm{NS}(A)$ be a nondegenerate line bundle on an
abelian variety $A/\C $. By Mumford's Vanishing Theorem, there is
a unique integer $i(\cL )$ such that $H^{i(\cL )}(A,\cL )\not =0$
and $H^j(A,\cL )=0$ for all $j\not =i(\cL )$ (cf.\,\cite{Mu},
Chapter III, Section 16, p.\,150). The so-called index $i(\cL )$
does only depend on the class of $\cL $ in $\mathrm{NS}(A)$ and we
have $0\leq i(\cL )\leq g=\dim (A)$.

If $H$ is the hermitian form associated to a line bundle $\cL $,
the index $i(\cL )$ agrees with the number of negative eigenvalues
of $H$ (\cite{Mu}, Chapter III, Section 16, p.\,162). By the
Riemann-Roch Theorem, $|K(\cL )| = |\mbox{Ker }\varphi _{\cL
}:A\rightarrow \hat A| = \mbox{dim}(H^{i(\cL )}(A,\cL ))$. In
particular, $\cL $ is {\em principal} when $\mbox{dim}(H^{i(\cL
)}(A,\cL ))=1$. Finally, $\cL $ is a polarization, i.\,e., an
ample line bundle, if and only if $i(\cL )=0$.

Let $A$ have quaternionic multiplication by an hereditary order
$\cO $ in $B$. It is our aim to compute the index of a line bundle
$\cL $ on $A$ in terms of the quaternion $\mu =c_1(\cL )$. For
these purposes, we introduce the following notation. Let $\mu \in
B_0$ be a pure quaternion. It satisfies $\mu ^2+\delta =0$ for
some $\delta \in F^*$. For any immersion $\sigma :F\hookrightarrow
\R$, let $\nu _{\sigma }= \nu _{\sigma }(\mu )\in \mathrm{GL}
_2(\R )$ be such that $\nu _{\sigma } \mu ^{\sigma } \nu _{\sigma
}^{-1}=\omega _{\sigma }$, where

$$
\omega _{\sigma } =
\begin{cases}
    \begin{pmatrix} _{0}& _{\sqrt {\sigma (\delta )}}\\
_{-\sqrt {\sigma (\delta )}}&_{0}
\end{pmatrix} & \text{if } \sigma (\delta )>0, \\
    \begin{pmatrix} _{\sqrt {\sigma (-\delta )}}&_{0}\\_{0}&_{-\sqrt {\sigma (-\delta
)}}
\end{pmatrix} & \text{otherwise}.
  \end{cases}
$$

We say that $\mu $ has {\em positive} or {\em negative orientation
at $\sigma $} according to the sign of the real number
$\mathrm{det }(\nu _{\sigma })$. Although $\nu _{\sigma }$ is not
uniquely determined by $\mu $, the sign sgn$(\mathrm{det }(\nu
_{\sigma }))$ is. Thus, to any pure quaternion $\mu \in B_0$, we
can attach a signature

$$
\mathrm{sgn } (\mu )= \mathrm{sgn}(\mathrm{det } (\nu _{\sigma
}))\in \{ \pm 1\} ^n.
$$

Motivated by the following theorem, we say that a pure quaternion
$\mu $ is {\em ample} with respect to $A$ if it has the same
orientation as $v_0\in V=\Lie (A)$: $\mathrm{sgn }(\mu )=
\mathrm{sgn }(\mathrm{Im }(\tau _i))$. For any real immersion
$\sigma :F\hookrightarrow \R $, define the local archimedean index
$i_{\sigma }(\mu )$ of $\mu $ by

$$
i_{\sigma }(\mu ) =
  \begin{cases}
    0& \text{if } \sigma (\delta )>0
\mbox{ and det}(\nu _{\sigma })\cdot \mathrm{Im }(\tau _{\sigma
})>0,
\\
    1 & \text{if } \sigma (\delta )<0, \\
    2 & \text{if } \sigma (\delta )>0
\mbox{ and det}(\nu _{\sigma })\cdot \mathrm{Im }(\tau _{\sigma
})<0.
  \end{cases}
$$

With these notations we have

\begin{theorem}
\label{index}

The index of $\cL \in \mathrm{NS}(A)$ is
$$
i(\cL ) = \sum _{\sigma :F\hookrightarrow \R } i_{\sigma }(\mu ).
$$

\end{theorem}

{\em Proof. } The index of $i(\cL )$ coincides with the number of
negative eigenvalues of the hermitian form $H_{\mu }$ associated
to the line bundle $\cL $. If we regard $\M _2(\R )\times
\stackrel {n}{...}\times \M _2(\R )$ embedded diagonally in $\M
_{2n}(\R )$, there is an isomorphism of real vector spaces between
$B\otimes _{\Q }\R $ and $\M _2(\R )\times \stackrel
{n}{...}\times \M _2(\R )$ explicitly given by the map $\beta
\mapsto \beta \cdot v_0$. The complex structure that $\M _2(\R
)^n$ inherits from that of $V$ is such that $\{ 0\} \times
...\times \M _2(\R )\times ...\times \{ 0\} $ are complex vector
subspaces of $\M _2(\R )^n$ and we may choose a $\C $-basis of $V$
of the form $\{ \mbox{diag }(\beta _1, 0, ..., 0)\cdot v_0, $
$\mbox{diag }(\gamma _1, 0, ..., 0)\cdot v_0, ..., \mbox{diag }(0,
..., 0,\beta _n)\cdot v_0,$ $\mbox{diag }(0, ..., \gamma _n)\cdot
v_0\} $ for $\beta _i$, $\gamma _i\in \M _2(\R )$.

Let $\iota =\mathrm{diag} (\iota _{\sigma })\in \mathrm{GL}
_{2n}(\R )$ be such that $\iota \cdot v_0=\sqrt{-1} v_0$. For any
$\beta =\mbox{diag }_{\sigma } (\beta _{\sigma })$, we have that
$\gamma =\mbox{diag }_{\sigma } (\gamma _{\sigma })\in \M _{2n}(\R
)$ and

$$
H_{\mu }(\beta v_0, \gamma v_0)= \sum _{\sigma }\mbox{tr }(\mu
^{\sigma }\beta _{\sigma } \iota _{\sigma }\bar \gamma _{\sigma })
+ \sqrt{-1}\sum _{\sigma } \mbox{tr }(\mu ^{\sigma }\beta _{\sigma
} \bar {\gamma _{\sigma }}).
$$

Thus, if we let $H_{\sigma }\in \mathrm {GL}_2(\C )$ denote the
restriction of $H_{\mu }$ to $V_{\sigma }=\M _2(\R )\cdot
{\begin{pmatrix} {\tau _{\sigma }} & {1} \end{pmatrix}}^t$, the
matrix of $H_{\mu }$ respect to the chosen basis has diagonal form
$H_{\mu }=\mbox{diag }(H_{\sigma })$.

In order to prove Theorem \ref{index}, it suffices to show that
the hermitian form $H_{\sigma }$ has $i_{\sigma }(\mu )$ negative
eigenvalues. Take $\beta \in \M _2(\R )$ and let $v=\beta \cdot
{\begin{pmatrix} {\tau _{\sigma }} & {1}\end{pmatrix}}^t\in
V_{\sigma }$. Then, $H_{\sigma }(v, v)= \mathrm{tr }(\mu ^{\sigma
}\beta _{\sigma }\iota _{\sigma }\bar \beta _{\sigma })=\mathrm{tr
}(\omega _{\sigma }\beta '_{\sigma } \iota '_{\sigma }\bar \beta
'_{\sigma } )$, where $\beta '_{\sigma }=\nu _{\sigma }\beta
_{\sigma } \nu _{\sigma }^{-1}$ and $\iota '_{\sigma }=\nu
_{\sigma }\iota _{\sigma }\nu _ {\sigma }^{-1}$. Denote $w
_{\sigma }= {\begin{pmatrix} w_1 & w_2 \end{pmatrix}}^t = \nu
_{\sigma } \beta _{\sigma }\cdot {\begin{pmatrix} \tau _{\sigma }
& 1
\end{pmatrix}}^t \in \C ^2$ and $\| w_{\sigma }\| ^2 = w_1 \bar w_1+w_2 \bar w_2$.

Some computation yields that

$$
H_{\sigma }(v, v)=\sum _{\sigma } \dfrac{C_{\sigma }\sqrt {|
\sigma (\delta )|}}{\mathrm{det}(\nu _{\sigma }) \mathrm{Im} (\tau
_{\sigma })},
$$
where $C_{\sigma }=\| w\| ^2 $ if $\sigma (\delta )>0$ and
$C_{\sigma }= 2\mathrm{ Re}(w_1 \bar w_2)$ if $\sigma (\delta
)<0$. From this, the result follows. $\Box $

\begin{remark}
From the above formula, the well-known relation $ i(\cL )+i(\cL
^{-1})=\dim A $ (\cite{Mu}, Chapter III, Section 16, p.\,150) is
reobtained.
\end{remark}

\section{{\bfseries Principal polarizations and self-duality}}

If an abelian variety $A$ admits a principal line bundle, and
hence is self-dual, it is natural to ask whether it admits a
principal polarization. It is the purpose of this section to study
this question under the assumption that $A$ has quaternionic
multiplication by an hereditary order $\cO $.

From Corollary \ref{dual}, a sufficient condition for $A$ to be
self-dual when $(\vartheta _{F/\Q }, \cD _{\cO }) = 1$ is that
$\cD _{\cO }$ and $\mathrm{n}(\cI )\cdot \vartheta _{F/\Q }$ are
principal ideals. By Theorem \ref{index}, a necessary condition
for $A$ to be principally polarizable is that $\cD _{\cO }$ be
generated by a totally positive element $D$ in $F$. However, in
general this is not sufficient for the existence of a principal
polarization on $A$.

Let $F^*_+$ denote the subgroup of totally positive elements of
$F^*$, $R_{F_+}^* = R_F^*\cap F^*_+$ and $\cO ^*_+ = \{ \alpha \in
\cO ^*: \n (\alpha )\in R_{F_+}^* \} $. Let $\Pic _+(F)$ be the
narrow class group of $F$ and let $h_+(F)=|\Pic _+(F)|$. We let
$\Sigma = \Sigma (R_F^*)\subseteq \{ \pm 1\} ^n$ be the $\F
_2$-subspace of signatures of units in $R_F^*$. As $\F _2$-vector
spaces, $\Sigma \simeq R_F^*/R_{F +}^*$ and, by Dirichlet's Unit
Theorem, $|\Sigma | = \dfrac {2^n h(F)}{h_+(F)}$. Since $\n (\cO
^*)=R_F^*$, the group $\Sigma $ fits in the exact sequence

$$
1\rightarrow \cO ^*_+\rightarrow \cO ^*\rightarrow \Sigma
\rightarrow 1
$$
$$
\qquad \qquad \qquad \qquad \qquad \alpha \mapsto (\mathrm{sgn}
(\mathrm{det } \alpha ^{\sigma })).\qquad
$$

\begin{definition}

We denote by $\Omega \subseteq \{\pm 1\} ^n$ the set of signatures
$$
\Omega = \{ (\mathrm{sgn }(\mathrm{det } \nu _{\sigma }(\mu )) :
\mu \in \cP (\cO ) \} .
$$

\end{definition}

The set $\Omega $ can be identified with the set of connected
components of $\R ^n\setminus \cup _{i=1}^n\{x_i = 0\}$. With the
notations as above, we obtain the following corollary of Theorems
\ref{pral} and \ref{index}.

\begin{proposition}
\label{polpral}

Let $\cI $ be an ideal of $B$ and assume that its left order $\cO
$ is hereditary. Then, there exist principally polarizable abelian
varieties $A$ with quaternionic multiplication by $\cO $ and
$H_1(A,\Z )\cong \cI $ if and only if $\cD _{\cO }$ and
$\mathrm{n}(\cI )\cdot \vartheta _{F/\Q }$ are principal ideals
and $\cD _{\cO }=(D)$ is generated by an element $D\in F_+^*$.

If this is the case, an abelian variety $A=V/\cI \cdot (\tau _1,
1, ..., \tau _n, 1)^t$, admits a principal polarization if and
only if $(\mathrm{sgn}(\mathrm{Im }\tau _i))\in \Omega $.

\end{proposition}

We note in passing that, as a consequence of the above corollary
and an application of \v{C}ebotarev's Density Theorem (\cite{Ne},
Chapter VII, Section 13), self-dual but non principally
polarizable abelian varieties $A$ can be constructed. Examples of
these abelian varieties are nontrivial, since in the generic case
in which the ring of endomorphisms is $\End (A)=\Z $, the abelian
variety $A$ is principally polarizable if and only if it is
self-dual.

Signature questions on number fields are delicate. In order to
have a better understanding of Proposition \ref{polpral}, we
describe $\Omega $ as the union (as sets) of linear varieties in
the affine space $\mathbb A _{\F _2}^n = \{ \pm 1\} ^n$ as
follows. Let $\{ u_k\} $ be a set of representatives of units in
$R_F^*/R_F^{*2}$ and, for any order $S\supseteq R_F[\sqrt {-u_k
D}]$ in $L=F(\sqrt {-u_k D})$, choose $\mu _S\in \cP (S, \cO )$.
We considered in Section \ref{poll} the Galois group $\Delta =
\mbox{Ker}(N:\Pic(S)\rightarrow \Pic (F))$. Naturally associated
to it there is a sub-space of signatures $\Sigma (\Delta )$ in the
quotient space $\mathbb A_{\F _2}^n/\Sigma (R_F^*)$: if $\mathfrak
b$ is an ideal of $S$ such that $N_{L/F}(\mathfrak b)=(b)$ for
some $b\in F^*$, the signature of $b$ does not depend on the
choice of $\mathfrak b$ in its class in $\Pic (S)$ but depends on
the choice of the generator $b$ up to signatures in $\Sigma
(R_F^*)$. By an abuse of notation, we still denote by $\Sigma
(\Delta )$ the sub-space of $\mathbb A^n_{\F _2}$ generated by
$\Sigma (R_F^*)$ and the signatures of the norms of ideals in
$\Delta $. Then, from Proposition \ref{Galois} we obtain that

$$
\Omega = \sqcup _{k, S} \Sigma (\Delta )\cdot \mathrm{sgn } (\mu
_S),
$$
as a disjoint union. This allows us to compute the set $\Omega $
in many explicit examples and to show that, in many cases,
coincides with the whole space of signatures $\{ \pm 1\} ^n$. The
following corollary, which remains valid even if we remove the
assumption $(\vartheta _{F/\Q }, \cD _{\cO }) = 1$, illustrates
this fact.

\begin{corollary}
\label{cor}

Let $F$ be a totally real number field of degree $[F:\Q ]=n$. Let
$\cO $ be an hereditary order in a totally indefinite quaternion
algebra $B$ over $F$ and let $\cI $ be a left $\cO $-ideal such
that $\cD _{\cO } = (D)$ for $D\in F^*_+$ and $\n (\cI )\cdot
\vartheta _{F/\Q } = (t^{-1})$ for $t\in F^*$.

If the narrow class number $h_+(F)$ of $F$ equals the usual class
number $h(F)$, i.\,e. if $\Sigma (R_F^*)=\{ \pm 1\} ^n$, then any
abelian variety $A$ with quaternionic multiplication by $\cO $ and
$H_1(A, \Z )\cong \cI $ is principally polarizable.

In particular, if $h_+(F)=1$, then the above conditions on $\cO $
and $\cI $ are accomplished.

\end{corollary}

{\em Proof. } Since $\Sigma (R_F^*) = \{ \pm 1\} ^n$, we have
$\Omega =\{ \pm 1\} ^{n}$ and the result follows from Proposition
\ref{polpral}. $\Box $

This is highly relevant in the study of certain {\em Shimura
varieties}. As was already known to the specialists in the case of
maximal orders, we obtain that {\em any} abelian surface with
quaternionic multiplication by an {\em hereditary} order in an
indefinite quaternion algebra $B/\Q $ admits a principal
polarization.

\section{{\bfseries The number of isomorphism classes of principal polarizations}}

Let $A=V/\Lambda $ with $\Lambda = \cI \cdot v_0$ be an abelian
variety with quaternionic multiplication by a maximal order $\cO
$. For any integer $0\leq i\leq g$, let $\Pi _i(A)$ denote the set
of isomorphism classes of principal line bundles $\cL \in
\mathrm{NS}(A)$ of index $i(\cL )=i$. The set $\Pi (A)$ naturally
splits as the disjoint union $\Pi (A)= \sqcup \Pi _i(A)$.
Moreover, due to the relation $i(\cL ) + i(\cL ^{-1}) = g$, the
map $\cL \mapsto \cL ^{-1}$ induces a one-to-one correspondence
between $\Pi _i(A)$ and $\Pi _{g-i}(A)$.

Formulas for $\pi _i(A)$, $0\leq i\leq g$, analogous to that of
Theorem \ref{pi} can be derived. Due to its significance, we will
only concentrate on the number $\pi _0(A)$ of classes of principal
polarizations. The Galois action on the sets $E (S, \cO )$ of
Eichler classes of optimal embeddings and its behaviour respect to
the index of the associated line bundles will play an important
role.

Assume then that $\Pi _0(A)\not =\emptyset $. For simplicity,
recall that we also assume that $(\vartheta _{F/\Q }, \cD _{\cO
})=1$. By Proposition \ref{polpral}, we can choose $D\in F^*_+$
and $t \in F^*$ such that $\cD _{\cO }=(D)$ and $\mathrm{n}(\cI
)\cdot \vartheta _{F/\Q }=(t ^{-1})$.

Let $u\in R_{F_+}^*$ be a totally positive unit. Let us agree to
say that an order $S\supseteq R_F[\sqrt {-u D}]$ is {\em ample}
respect to $\cO $ if there exists an optimal embedding
$i:S\hookrightarrow \cO $ such that $\mu =i(\sqrt {-u D})$ is
ample (cf.\,the discussion preceding Theorem \ref{index}). Define
$\cS _u$ to be the set of ample orders $S\supseteq R_F[\sqrt {-u
D}]$ in $F(\sqrt {-u D})$. The existence of a principal
polarization $\cL $ on $A$ implies that there is some $\cS _u$
nonempty. With this notation, we obtain the following expression
for $\pi _0(A)$ in terms of the narrow class number of $F$ and the
class numbers of certain CM-fields that embed in $B$.

\begin{theorem}
\label{pio}

The number of isomorphism classes of principal polarizations on
$A$ is

$$
\pi _0(A) = \frac {1}{2 h_+(F)}\sum _{u\in R_{F_+}/R_F^{*2}} \sum
_{S\in \cS _u}
 2^{e^+_S} h(S),
$$
where $2^{e_S^+} = |R_{F +}^*/N(S^*)|$.
\end{theorem}

{\em Proof. } By the existing duality between $\Pi _0(A)$ and $\Pi
_g(A)$, it is equivalent to show that $\pi _0(A) + \pi _g(A) =
\sum _u \sum _{S\in \cS _u} 2^{e^+_S} h(S)/h_+(F)$.

Let us introduce the set $\cP _{0, g}(\cO )  = \{ \mu \in \cO :
\mathrm{sgn } (\mu ) = \pm (\mathrm{sign } (\mathrm{Im } \tau
_i)), \n (\mu )\in R_{F_+}^*\cdot D \} $. By Theorems \ref{NeSe}
and \ref{index}, the set $P _{0, g}(\cO )= \cP _{0, g}(\cO ) ^{\pm
}/_{\sim _p}$ is in one-to-one correspondence with $\Pi_0(A)\cup
\Pi _g(A)$ and we have a natural decomposition $P _{0, g}(\cO ) =
\sqcup P _{0, g}(S, \cO )$ as $S$ runs among ample orders in $\cS
_u$ and $u\in R_{F_+}^*/R_F^{*2}$.

Fix $u\in R^*_{F_+}$ and $S$ in $\cS _u$. In order to compute the
cardinality of $P _{0, g}(S, \cO )$, we relate it to the set
$E_{0, g}(S, \cO )=\cP _{0, g}(S, \cO )/_{\sim _{e}}$ of $\cO
^*_{\pm }$-Eichler conjugacy classes of optimal embeddings $i_{\mu
}:S\hookrightarrow \cO $. Here, we agree to say that two
quaternions $\mu _1$ and $\mu _2\in \cP _{0, g}(S, \cO )$ are
Eichler conjugated by $\cO ^*_{\pm }$ if there is a unit $\alpha
\in \cO ^*_{\pm }$ of either totally positive or totally negative
reduced norm such that $\mu _1 = \alpha ^{-1} \mu _2\alpha $. Note
that, by Theorem \ref{index}, the action of $\cO ^*_{\pm
}$-conjugation on the line bundle $\cL $ associated to an element
in $\cP (S, \cO ) $ either preserves the index $i(\cL )$ or
switches it to $g-i(\cL )$. This makes sense of the quotient $E
_{0, g}(S, \cO )$.

We have the following exact diagram:

$$
\begin{matrix}

0   &\ra & \Delta _+   & \ra   & \Pic (S)   & \stackrel
{N_{L/F}}{\ra } & \Pic
_+(F)  & \ra  & 0  \\
      & & \downarrow   & & \parallel   & & \downarrow
      & & \\

 0    & \ra  & \Delta    & \ra   & \Pic (S)   & \stackrel {N_{L/F}}{\ra } &
\Pic (F)  & \ra  & 0.

\end{matrix}
$$

Indeed, there is a natural map $\Pic (S)\stackrel
{N_{L/F}}{\rightarrow } \Pic _+(F)$, since the norm of an element
$a+b\sqrt {-u D}\in L$ for $a$, $b\in F$ is $a^2+u b^2D\in F^*_+$.
The surjectivity of the map $\Pic (S)\rightarrow \Pic _+(F)$ is
argued as in Section \ref{poll} by replacing the Hilbert class
field $H_F$ of $F$ by the big Hilbert class field $H^+_F$, whose
Galois group over $F$ is $\mathrm{Gal} (H_F^+/F)=\Pic _+(F)$. By
Proposition \ref{Galois}, $\Delta $ acts freely and transitively
on $E (S, \cO )$. Therefore, by Theorem \ref{index}, there is also
a free action of $\Delta _+$ on $E _{0, g}(S, \cO )$. Up to sign,
the $\cO ^*_{\pm }$-Eichler conjugation class of an element $\mu
\in \cP (S, \cO )$ has a well-defined orientation $\pm \mathrm{sgn
} (\mu)$. Note also that two inequivalent $\cO ^*_{\pm }$-Eichler
classes that fall in the same $\cO ^*$-conjugation class are never
oriented in the same manner, even not up to sign. Taken together,
this shows that $\Delta _+$ also acts transitively on $E _{0,
g}(S, \cO )$. This means that

$$
|E _{0, g}(S, \cO )| = \frac {h(S)}{h_+(F)}.
$$

There is again a natural surjective map $\rho :P _{0, g}(S, \cO
)\ra E _{0, g}(S, \cO )$ and, arguing as in Section \ref{proofpi},
Theorem \ref{pio} follows. $\Box $

\vspace{0.5cm} {\bfseries Examples in low dimensions. } When we
particularize our results to dimension $2$, we obtain an easy to
apply expression for the number of principal polarizations on an
abelian surface with maximal quaternionic multiplication, as it is
stated in Theorem \ref{main} in the introduction. As an example,
the number of isomorphism classes of principal polarizations on an
abelian surface $A$ with quaternionic multiplication by a maximal
order in a quaternion algebra of discriminant $D = 2\cdot 3\cdot
5\cdot 7\cdot 11\cdot 13\cdot 17\cdot 19$ is $\pi _0(A)=1040$.
This also implies the existence of $1040$ pairwise nonisomorphic
smooth algebraic curves $C_1, ..., C_{1040}$ of genus $2$ such
that their Jacobian varieties are isomorphic as unpolarized
abelian surfaces.

In addition, since $\pi (A)=\pi _0(A)+\pi _1(A)+\pi _2(A)$ and
$\pi _0(A)=\pi _2(A)$, Theorems \ref{pi} and \ref{pio} yield the
formula

$$
\pi _1(A) =
  \begin{cases}
    \varepsilon _{4 D}\cdot h(4 D)+ \varepsilon _D\cdot h(D) &
\text{ if }D\equiv 1\quad \mbox{ mod }4, \\
    \varepsilon _{4 D}\cdot h(4 D) & \text{otherwise,}
  \end{cases}
$$
where $\varepsilon _{D}$ and $\varepsilon _{4 D} = 1$ or $\frac
{1}{2}$ is computed from the formula for $e_S$ in Theorem
\ref{pi}.

Let $F$ be the real quadratic field $\Q (\sqrt {2})$ and let $B$
be the quaternion algebra over $F$ that ramifies exactly at the
two prime ideals $(3\pm \sqrt {2})$ above $7$. By applying
Theorems \ref{pi}, \ref{index} and \ref{pio} and the valuable help
of the programming package PARI (\cite{PA}), we conclude that, for
any abelian four-fold $A$ such that $\End (A)$ is a maximal order
in the quaternion algebra $B/\Q (\sqrt {2})$ of discriminant $7$,
the number of isomorphism classes of principal line bundles of
index 0, 1, 2, 3 and 4 are $\pi _0(A)=\pi _4(A)=6$, $\pi _1(A)=\pi
_3(A)=4$ and $\pi _2(A)=4$, respectively.

\section{{\bfseries Asymptotic behaviour of $\pi _0(A)$}}\label{Asymp}

We can combine Theorem \ref{pio} with analytical tools to estimate
the asymptotic behaviour of $\log (\pi _0(A))$. This will yield a
stronger version of part 1 of Theorem \ref{ass} in the
introduction. For any number field $L$, we let $D_L$ and
$\mathrm{Reg }_L$ stand for the absolute value of the discriminant
and the regulator, respectively.

\begin{theorem}
\label{log}

Let $F$ be a totally real number field of degree $n$. Let $A$
range over a sequence of principally polarizable abelian varieties
with quaternionic multiplication by a maximal order in a totally
indefinite quaternion algebra $B$ over $F$ of discriminant $D\in
F^*_+$ with $|\mbox{N}_{F/\Q }(D)|\rightarrow \infty $. Then

$$
\log \pi _0(A)\sim \log \sqrt {|\mbox{N}_{F/\Q }(D)|\cdot D_F}.
$$

\end{theorem}

The proof of Theorem \ref{log} adapts an argument of Horie-Horie
(\cite{HoHo}) on estimates of relative class numbers of CM-fields.
We first show that there indeed exist families of abelian
varieties satisfying the properties quoted in the theorem.

By \v{C}ebotarev's Density Theorem, we can find infinitely many
pairwise different totally positive principal prime ideals $\{\wp
_i\}_{i\geq 1}$ in $F$. We can also choose them such that $(\wp
_i, \vartheta _{F/\Q })=1$. We then obtain principal ideals $(D_j)
= \wp _1\cdot \wp _2\cdot ...\cdot \wp _{2 j - 1}\cdot \wp _{2 j}$
with $D_j\in F^*_+$ and $(D_j, \vartheta _{F/\Q })=1$. According
to \cite{Vi}, p.\,74, there exists a totally indefinite quaternion
algebra $B_j$ over $F$ of discriminant $D_j$ for any $j\geq 1$.
Then, Proposition \ref{polpral} asserts that there exists an
abelian variety $A_j$ of dimension $2n$ such that $\End (A_j)$ is
a maximal order in $B_j$ and $\Pi _0(A_j)\not = \emptyset $.

$\\ ${\em Proof of Theorem \ref{log}.} Let $A$ be a principally
polarizable abelian variety with quaternionic multiplication by a
maximal order in a totally indefinite division quaternion algebra
$B$ over $F$ of discriminant $D\in F^*_+$.

For any totally positive unit $u_k\in R_{F_+}^*$, let $L_k=
F(\sqrt {-u_k D})$. For any order $S\supseteq R_F[\sqrt {-u_k
D_j}]$ in the CM-field $L_k$, it holds that $h(S) = c_S h(L_k)$
for some positive constant $c_S\in \Z $ which is uniformly bounded
by $2^n$.

The class number $h(F)$ turns out to divide $h(L_k)$ and the {\em
relative class number} of $L_k$ is defined to be $h^-(L_k) =
h(L_k)/h(F)$ (cf.\,\cite{Lo1}). Since $h_+(F) = 2^m h(F)$ for $m =
n - \mathrm{dim }_{\F _2}(\Sigma (R_F^*))$, Theorem \ref{pio} can
be rephrased as $\pi _0(A) = \sum 2 ^{(e^+_S - 1 -m)} c_S
h^-(L_k)$.

In order to apply the Brauer-Siegel Theorem, the key point is to
relate the several absolute discriminants $D_{L_k}$ and regulators
$\mathrm{Reg } _{L_k}$ as $u_k$ vary among totally positive units
in $F$. Firstly, we have the relations $D_{L_k} = |N_{F/\Q
}(D_{{L_k}/F})\cdot D_F^2| = 2^{p_k} |N_{F/\Q }(D)| D_F^2$ for
some $0\leq p_k \leq 2 n$. Secondly, by \cite{Wa}, p.\,41, it
holds that $\mathrm {Reg }_{L_k}=2^{c} \mathrm {Reg }_F$ with
$c=n-1$ or $n-2$.

Let $\varepsilon $ be a sufficiently small positive number. By the
Brauer-Siegel Theorem, it holds that $D_{L_k}^{(1-\varepsilon
)/2}\leq h(L_k)\mathrm {Reg }_{L_k}\leq D_{L_k}^{(1+\varepsilon
)/2}$ for $D_{L_k}\gg 1$. Thus
$$
\frac{D_F^{(1-\varepsilon )/2}}{h(F)\mathrm {Reg }_F}
(D_{L_k}/D_F) ^{(1-\varepsilon )/2}\leq h^-(L_k)\leq
\frac{D_F^{(1+\varepsilon )/2}}{h(F)\mathrm {Reg }_F}
(D_{L_k}/D_F) ^{(1+\varepsilon )/2}.
$$
Fixing an arbitrary CM-field $L$ in the expression for $\pi
_0(A)$, this boils down to
$$
C_-\cdot \frac{D_F^{(1-\varepsilon )/2}}{h(F)\mathrm {Reg }_F}
(D_{L}/D_F) ^{(1-\varepsilon )/2}\leq \pi _0(A)\leq C_+\cdot
\frac{D_F^{(1+\varepsilon )/2}}{h(F)\mathrm {Reg }_F} (D_{L}/D_F)
^{(1+\varepsilon )/2}
$$
for some positive constants $C_-$ and $C_+$. Taking logarithms,
these inequalities yield Theorem \ref{log}. $\Box $

\begin{remark}

The argument above is not effective since it relies on the
classical Brauer-Siegel Theorem on class numbers. However, recent
work of Louboutin (\cite{Lo1}, \cite{Lo2}) on lower and upper
bounds for relative class numbers of CM-fields, based upon
estimates of residues at $s=1$ of Dedekind zeta functions, could
be used to obtain explicit lower and upper bounds for $\pi _0(A)$.
\end{remark}

Finally, we conclude this paper with the proof of the second main
result quoted in the Introduction.

$\\ ${\em Proof of Theorem \ref{ass}. } Part $1$ is an immediate
consequence of Theorem \ref{log}. Let us explain how part $2$
follows. Assume that $A$ is a simple complex abelian variety of
odd and square-free dimension $g$. Then, by Albert's
classification of simple division algebras (\cite{Mu}, Chapter IV,
Section 19 and 21), $\End (A)\simeq S$ is an order in either a
totally real number field $F$ or a CM-field $L$ over a totally
real number field $F$. In any case, $[F:\Q ]\leq g$. In the former
case, by Theorem 3.1 of Lange in \cite{La}, $\pi
_0(A)=|S^*_+/S^{*2}|\leq 2^{g-1}$. In the latter, let $S_0\subset
F$ be the subring of $S$ fixed by complex conjugation. If $\cL $
is a principal polarization on $A$, the Rosati involution
precisely induces complex conjugation on $\End (A)\simeq S$ and we
have that $\pi _0(A)= |S_{0 +}^*/\Norm _{L/F}(S^*)|\leq |S_{0
+}^*/S_0^{*2}|\leq 2^{g-1}$, by applying Theorem 1.5 of \cite{La}.

\end{document}